\documentclass[11pt]{article}
\usepackage{amsthm,epsfig}
\def\E{{\rm Exp}}
\def\P{{\rm Pr}}
\def\proof{\smallskip\noindent{\bf Proof:} }

\newtheorem{theorem}{Theorem}
\newtheorem{corollary}{Corollary}[theorem]
\newtheorem{lemma}{Lemma}[section]

\begin{document}

\title{Construction of locally plane graphs with many edges}
\author{
G\'abor Tardos\footnote{Reserach partially supported by NSERC grant 329527 and
  by OTKA grants T-046234, AT048826 and NK-62321.}\\
R\'enyi Institute, Budapest, Hungary and\\
School of Computing Science, Simon Fraser University, Burnaby, BC\\
{\tt tardos@renyi.hu}
}
\date{}

\maketitle

\begin{abstract}
A graph drawn in the plane with straight-line edges is called a
geometric graph. If no path of length at most $k$ in a geometric graph $G$ is
self-intersecting we call $G$ $k$-locally plane. The main result of
this paper is a construction of $k$-locally plane graphs with a
super-linear number of edges. For the proof we develop randomized thinning
procedures for edge-colored bipartite (abstract) graphs that can be applied to
other problems as well.
\end{abstract}

\section{Introduction}

A {\em geometric graph} $G$ is a straight-line drawing of a simple, finite
(abstract) graph $(V,E)$, i.e., we identify the vertices $x\in V$ with
distinct points in the Euclidean plane, and we identify any edge $\{x,y\}\in E$
with the straight line segments $xy$ in the plane. We assume that the edge
$xy$ does not pass through any vertex of $G$ besides $x$ and $y$. We call
$(V,E)$ the abstract graph {\em underlying} $G$. We say
that the edges $e_1,e_2\in E$ {\em cross} if the corresponding line
segments cross each other, i.e., if they have a common interior point. We
say that a subgraph of $G$ is {\em self-intersecting} if it contains a
pair of crossing edges.

Geometric graphs without crossing edges are plane drawings of planar
graphs: they have at most $3n-6$ edges if $n\ge3$ is the number of vertices.

Avital and Hanani \cite{AH}, Erd\H os, and Perles initiated in the mid 1960s
the systematic study of similar questions for more complicated {\em forbidden
configurations}: Let $H$ be set of
forbidden configurations (geometric subgraphs). What is the maximal number of
edges of an $n$ vertex geometric graph not containing any configuration
belonging to $H$? This problem can be regarded as a geometric version of the
fundamental problem of extremal graph theory: What is the maximum number of
edges that an abstract graph on $n$ vertices can have without containing
subgraphs of a certain kind.

Many questions of the above type on geometric graphs have been addressed in
recent years. In a number of papers linear upper bounds have been established
for the number of edges of a graph, under various forbidden configurations.
They include the configurations of three pairwise crossing edges \cite{AAPPS},
four pairwise crossing edges \cite{Ack},
the configurations of an edge crossed by many edges \cite{PT}, or even
two large stars with all edges of one of them crossing all edges of
the other \cite{TT}.

For a constant number of 5 or more pairwise crossing
edges Pavel Valtr has the best result \cite{Valtr}: a geometric graph
on $n$ vertices avoiding this configuration has $O(n\log n)$ edges but no
construction is known with a super linear number of edges.
Adam Marcus and the present author
\cite{MT} building on an earlier result of Pinchasi and Radoi\v ci\'c
\cite{PR} prove an $O(n^{3/2}\log n)$ bound on the number of edges of an
$n$ vertex geometric graph not containing self-intersecting cycles of length
four. No
construction is known beating the $O(n^{3/2})$ edges an abstract graph
having no cycles of length four can have.

For surveys on geometric graph theory, consult \cite{P1}, \cite{P2} and
\cite{PRT}.

In this paper we consider forbidding self-intersecting paths.
For $k\ge3$ we call a geometric graph {\em $k$-locally plane} if it has no
self-intersecting subgraph (whose underlying abstract graph is) isomorphic to
a path of length at most $k$.

Pach et al.\ \cite{PPTT} consider $3$-locally plane graphs, i.e., the case of
geometric graphs with no self-intersecting paths of length three. They prove matching lower
and upper bounds of $\Theta(n\log n)$ on the maximal number of edges of a
$3$-locally plane graph on $n$ vertices.

We extend the the lower bound result of \cite{PPTT} by considering
self-intersecting drawings of longer paths as forbidden configurations. Technically $k$-locally plane
graphs are defined by forbidding self intersecting paths of length $k$ {\em or
shorter}, but forbidding only self-intersecting paths of length {\em exactly
$k$} would lead to almost the same extremal function. Indeed, one can delete
at most $nk$ edges from any graph on $n$ vertices, such that all the non-zero
degrees in the remaining graph are larger than $k$. This ensures that all
shorter
paths can be extended to a path of length $k$. It is possible, but not likely,
that if one only forbids paths of length $k$ with {\em the first and last
edges crossing} significantly higher number of edges is achievable.

For even $k$ a geometric graph is $k$-locally plane if and only if the
$k/2$-neighborhood of any vertex $x$ is intersection free. Note that this
requirement is much stronger than the similar condition on abstract
graphs, namely that the $k/2$ neighborhood of any point is planar.
One can construct graphs with girth larger than $k$ and $\Omega(n^{k\over
k-1})$ edges. In such a graph the $k/2$-neighborhood of any vertex is a tree,
still by \cite{PPTT} the graph does not even have $3$-locally plane drawing.

Extending the lower bound result in \cite{PPTT} we prove in Theorem~\ref{15}
that for arbitrary
fixed $k\ge3$ there exist $k$-locally plane
graphs on $n$ vertices with $\Omega(n\log^{(\lfloor
k/2\rfloor)}n)$ edges. Here $\log^{(t)}$ denotes $t$
times iterated logarithm and the hidden constant in $\Omega$
depends on $k$. Given two arbitrarily small
disks in the plane we can even ensure that all edges of the constructed graph
connect a vertex from the first disk with another vertex from the second. This
ensures that all the edges of the constructed geometric graph are
arbitrarily close to each other in length and direction. In the view of the
author this makes the existence of a high average degree (or for that matter
high minimum degree) 100-locally plane graphs even more surprising.

As a simple corollary we can characterize the abstract graphs $H$ such that
any geometric graph having no self-intersecting subgraph isomorphic to $H$ has
a linear number of edges. These graphs $H$ are the forests with at least two
nontrivial components. To see the linear bound for the number of edges of a
geometric graph avoiding a self-intersecting copy of such a forest $H$
first delete a linear number of edges from an arbitrary geometric graph $G$
until all non-zero degrees of the remaining geometric graph $G'$ are at least
$|V(H)|$. If $G'$ is crossing free the linear bound of the number of edges
follows. If you find a pair of crossing edges in $G'$ they can be extended to
a subgraph of $G'$ isomorphic to $H$. On the other hand, if $H$ contains a
cycle, then even an abstract graph avoiding it can have a super-linear
number of edges. If $H$ is a tree of diameter $k$, then a $k$-locally plane
geometric graph has no self-intersecting copy of $H$. Notice that
the extremal number of edges in this case (assuming $k>2$) is $O(n\log n)$ by
\cite{PPTT}, thus much smaller than the $\Omega(n^\alpha)$ edges ($\alpha>1$)
for forbidden cycles.

The main tool used in the proof of the above result is a randomized thinning
procedure that takes a $d$ edge colored bipartite graph of average degree
$\Theta(d)$ and returns a subgraph on the same vertex set with average
degree $\Theta(\log d)$ that does not have a special type of colored
path (walk) of length four. The procedure can be applied recursively to
obtain a subgraph avoiding longer paths of certain types. We believe
this thinning procedure to be of independent interest. In particular it can be
used to obtain optimal 0-1 matrix constructions for certain avoided submatrix
problems, see the exact statement in Section \ref{s4} and the details in
\cite{T}.

In Section \ref{s2} we define two thinning procedures for edge colored
bipartite graphs and prove their main properties. This is the most technical
part of the paper. While these procedures proved useful in other setting too
(and the author finds the involved combinatorics appealing) this entire
section can be skipped if one reads the definition of $k$-flat graphs (the two
paragraphs before Lemma~\ref{10}) and is
willing to accept Corollary~\ref{corol} at the end of the section (we also use
the simple observation in Lemma~\ref{11}). In fact, in order to understand
the main ideas behind the main result of this paper it is recommended to skip
Section~\ref{s2} on the first reading and to go straight to Section~\ref{s3}
where we use Corollary~\ref{corol} to construct locally plane graphs
with many edges. In Section \ref{s4} we comment on the optimality of the
thinning procedures and have some concluding remarks.

\section{Thinning}\label{s2}

In this section we state and prove combinatorial statements about
edge colored abstract graphs, i.e., we do not consider here geometric graphs
at all. The connection to locally plane geometric graphs will be made clear in
Section~\ref{s3}. 

A {\em bipartite graph} is a triple $G=(A,B,E)$ with disjoint vertex sets
$A$ and $B$ (called {\em sides}) and edge set $E\subseteq A\times B$.
In particular, all
graphs considered in this paper are {\em simple}, i.e., they do not have
multiple edges or loops. The edge connecting the vertices $x$ and $y$ of $G$
is denoted by $xy$ or $(x,y)$. The latter notation is only used if $x\in A$
and $y\in B$. By a {\em $d$-edge coloring} of a graph
we mean a mapping $\chi:E\to\{1,2,\ldots d\}$ such that adjacent edges receive
different colors. When we do not specify $d$ we call such coloring a {\em
proper edge coloring} but we always assume that the ``set of colors'' are
linearly ordered. The degree of
any vertex in $G$ is at most the number $d$ of colors, and our
results are interesting if the average degree is close to $d$. Unless stated
otherwise the subgraphs of an edge colored graph are considered with the
inherited edge coloring. Our goal is to obtain a
subgraph of $G$ with as many edges as possible without containing a
certain type of colored path or walk.

\subsection{Heavy paths}

A simple example of the above concept is the following. We call a
path $P=v_0v_1v_2v_3$ of length $3$ {\em heavy} if $v_0\in B$ and
the colors $c_1=\chi(v_0v_1)$, $c_2=\chi(v_1v_2)$, $c_3=\chi(v_2v_3)$,
satisfy $c_2<c_1\le c_3$. The next lemma describes a thinning
procedure that gets rid of heavy paths. Although we do not need
this lemma in our construction, we present it as a simple analogue of our
results for more complicated forbidden walks.

\begin{lemma}\label0
Let $G=(A,B,E)$ be a bipartite graph with a proper edge coloring
$\chi:E\to\{1,2,\ldots d\}$. Then there exists a subgraph
$G'=(A,B,E')$ of $G$ with $|E'|\ge|E|/(3\lceil\sqrt d\,\rceil)$
that does not contain heavy paths.
\end{lemma}

The constant $3$ in the lemma could be replaced by the base of the natural
logarithm. Notice that if $G$ had average degree $\Theta(d)$, then the
average degree of $G'$ is $\Omega(\sqrt d)$.

\proof
Let $t=\lceil\sqrt d\,\rceil$ and select a uniform random value
$i_y\in\{1,2,\ldots,t\}$ independently for each vertex $y\in B$. We
say that an edge $e\in E$ is of {\em class} $\lceil \chi(e)/t\rceil$. We
call an edge $e=(x,y)\in E$ {\em eligible} if its class is $i_y$. Let
the subgraph $G'=(A,B,E')$ consist of those eligible edges $e=(x,y)\in E$
for which there exists no other eligible edge $(x,y')\in E$ of the same
class. Note that the words ``class'' and ``eligible'' will be used in a
different meaning when defining the two thinning procedures in the next
subsection. 

By the construction, all edges incident to a vertex $x\in A$ have
different classes and all edges incident to a vertex $y\in B$ have the
same class. Let $e_1$, $e_2$ and $e_3$ form a path in $G'$ starting in
$B$. Then $e_2$ and $e_3$ are of the same class, while the class of $e_1$ is
different. For the colors $c_i=\chi(e_i)$ this rules out the order $c_2<c_1\le
c_3$. Thus, $G'$ does not contain a heavy path. Note that another order,
$c_3\le c_1<c_2$ is also impossible.

The number of edges in $G'$ depends on the random choices we made.
Any edge $(x,y)\in E$ is eligible with probability $1/t$, and this is
independent for all the edges incident to a vertex $x\in A$. As
edges of a fixed color form a matching, there are at most $t$ edges
of any given class incident to $x$. Thus, we have
$$\P[(x,y)\in E']\ge\frac{(1-1/t)^{t-1}}t>\frac1{3t}.$$
The expected number of edges in $G'$ is
$$\E[|E'|]\ge\frac{|E|}{3t}.$$
It is possible to choose the random variables $i_y$ so that
the size of $E'$ is at least as large as its expected value. This proves the
lemma. \qed

\subsection{Fast and slow walks}

Next we turn to more complicated forbidden subgraphs. For motivation we
mention that self-crossing paths of length $4$ in the 3-locally plane graphs
of \cite{PPTT} (considered with their natural edge coloring) are exactly the
{\em fast walks} (to be defined below). For technical reasons, it
will be more convenient to consider walks, i.e., to permit that a vertex is
visited more than once, but we will not allow {\em backtracking}, i.e.,
turning back on the same edge immediately after it was traversed. Thus, for us
a {\em walk} of length $k$ is a sequence
$v_0,v_1,\ldots,v_k$ of vertices in the graph such that $v_{i-1}v_i$
is an edge for $1\le i\le k$ and $v_{i-2}\ne v_i$ for $2\le i\le
k$. The {\em $\chi$-coloring} (or simply coloring) of this walk is the
sequence $(\chi(v_0v_1),\chi(v_1v_2),\ldots,\chi(v_{k-1}v_k))$ of the colors of
the edges of the walk. If $\chi$ is a proper edge coloring, then any two
consecutive elements of the coloring sequence are different.

We use $\log$ to denote the binary logarithm. We introduce the notation
$P(a,b)$ for two non-equal strings
$a,b\in\{0,1\}^t$ to denote the first position $i\in\{1,2,\ldots,t\}$, where
$a$ and $b$
differ. We consider the set $\{0,1\}^t$ to be ordered lexicographically,
i.e., for $a,b\in\{0,1\}^t$ we have $a<b$ if $a$ has $0$ in position $P(a,b)$
(and thus $b$ has $1$ there).

The following trivial observation is used often in this paper. We state it
here without a proof.

\begin{lemma}\label1 Let $t\ge1$ and let $a$, $b$ and $c$ be distinct binary
strings of length $t$ with $P(a,b)<P(a,c)$. We have
$P(b,c)=P(a,b)$. Furthermore $a>b$ implies $c>b$, and $a<b$ implies $c<b$.
\end{lemma}

A walk of length $4$ with coloring $(c_1,c_2,c_3,c_4)$ is called a
{\em fast walk} if $c_2<c_3<c_4\le c_1$. Note that a fast
walk may start in either class $A$ or $B$. We call a walk of
length $4$ a {\em slow walk} if it starts in the class $B$ and
its coloring $(c_1,c_2,c_3,c_4)$ satisfies $c_2<c_3<c_4$
and $c_2<c_1\le c_4$. Note that either the color $c_1$ or
$c_3$ can be larger in a slow walk, or they can be equal.

The two {\em thinning procedures} below find a random subgraph of an
edge-colored bipartite graph. One is designed to avoid slow walks, the other
is designed to avoid fast walks.

\smallskip
\noindent{\bf Lexicographic thinning} Let $G=(A,B,E)$ be a
bipartite graph and let $\chi:E\to\{1,\ldots,d\}$ be a proper edge
coloring with $d\ge2$. {\em Lexicographic thinning} is a randomized procedure
that produces a subset $E'\subseteq E$ of the edges and the
corresponding subgraph $G'=(A,B,E')$ of $G$ as follows:

Let $t=\left\lceil\log d\over2\right\rceil+1$. Let $H$ be
the set of triplets $(a,i,z)$, where $a\in\{0,1\}^t$, $i\in\{2,3,\ldots,t\}$,
$z\in\{1,2,3,\ldots2^i\}$, and $a$ has $0$ in position $i$. Straightforward
calculation gives that $|H|=2^{2t}-2^{t+1}\ge2d$.

We order $H$ lexicographically, i.e., $(a,i,z)<(b,j,s)$ if $a<b$, or $a=b$ and
$i<j$, or $(a,i)=(b,j)$ and $z<s$.

Consider the following random function $F:\{1,\ldots,d\}\to H$.
We select uniformly at random the value
$F(1)=(a,i,z)\in H$ with the property that the first bit of $a$ is
$0$. We make $F(2)$ to be the next element of $H$ larger than $F(1)$,
and in general $F(k)$ is the next element of $H$ larger than $F(k-1)$ for
$2\le k\le d$. As $|H|\ge2d$ and $F(1)$ is chosen from the first half of $H$,
this defines $F$. In what follows we simply identify the color $k$ with the
element $F(k)\in H$ without any reference to the function $F$.

We say that $(a,i,z)\in H$ and any edge with this color is of {\em class} $a$
and {\em type} $i$, while $z$ will play no role except in counting how many
values it can take.

We choose an independent uniform random value $a_x\in\{0,1\}^t$ for each
vertex $x\in A\cup B$.
Let $e=(x,y)\in E$ be an edge of class $a$ and type $i$. We say that $e$ is
{\em eligible} if $a=a_y<a_x$ and $P(a,a_x)=i$.
Let the subgraph $G'=(A,B,E')$ contain those edges $e\in E$
that are eligible but not adjacent to another eligible edge $e'$ of the same
type as $e$.

\smallskip
\noindent{\bf Reversed thinning} Let $G=(A,B,E)$ be a bipartite
graph and let $\chi:E\to\{1,\ldots,d\}$ be a proper edge coloring with
$d\ge2$. {\em Reversed thinning} is a randomized procedure that produces a
subset $E'\subseteq E$ of the edges and the corresponding subgraph
$G'=(A,B,E')$ of $G$.

Reversed thinning is almost identical to the lexicographic thinning, the
only difference is in the ordering of the set $H$. We define $t$ and
$H$ as in the case of lexicographic thinning. Recall that $H$ is
the set of triplets $(a,i,z)$ where $a\in\{0,1\}^t$,
$i\in\{2,3,\ldots,t\}$, $z\in\{1,2,\ldots2^i\}$, and $a$ has $0$ in
position $i$. We still order $\{0,1\}^t$ lexicographically, but now
we reverse the lexicographic order of $H$ in the middle term $i$. That
is, we have $(a,i,z)<(b,j,s)$ if $a<b$, or $a=b$ and $i>j$ or
$(a,i)=(b,j)$ and $z<s$.

We define the function $F:\{1,\ldots,d\}\to H$, the {\em types} and {\em
classes} of colors and edges, {\em eligible} edges and the subset $E'$ of
edges the same way as for the lexicographic thinning, but using this modified
ordering of $H$.
\smallskip

Note that we associated a type in $\{2,\ldots,t\}$ and a class in $\{0,1\}^t$
to each edge in either procedure and they satisfy that
\begin{itemize}
\item an edge with a smaller class has smaller color;
\item among edges of equal class an edge with smaller type has smaller
  color in the case of lexicographic thinning and it has larger color in the
  case of reversed thinning;
\item among the edges incident to a vertex at most $2^i$ have the
  same class $a$ and the same type $i$.
\end{itemize}
Proving most of the properties of the thinning procedures this is all we need
to know about how classes and types are associated to the edges and we could
use a deterministic scheme for $F$. But for Lemma~\ref5 we need that all types
are well represented and the randomization in the identification function $F$
(as well as the dummy first bit of the class) is introduced to ensure this on
the average. This randomization is not needed if one assumes all color classes
have roughly the same size.

The next lemmas state the basic properties of the thinning
procedure. Lemma \ref2 lists common properties of the two procedures, while
Lemmas~\ref3 and \ref4 state the result of the thinning satisfies its ``design
criteria'' avoiding slow or fast walks. Finally Lemma~\ref5 shows that enough
edges remain in the constructed subgraphs on average. Note that Lemmas~\ref3
and \ref4 are special cases of the more complex Lemmas~\ref7 and
\ref8 proved independently. We state and prove the simple cases separately for
clarity, but these proofs could be skipped.

\begin{lemma}\label2 Let $G=(A,B,E)$ be a bipartite graph with a
proper edge coloring $\chi:E\to\{1,\ldots,d\}$. If
$G'=(A,B,E')$ is the result of either the lexicographic or the
reversed thinning then we have

\begin{description}
\item[a)]Adjacent edges in $G'$ have distinct types.
\item[b)]If two edges of $G'$ meet in $B$, they have
the same class.
\item[c)]Suppose two distinct edges $e$ and $e'$ of $G'$ meet in $A$. Let
their classes and types be $a$, $a'$ and $i$, $i'$, respectively. If $i<i'$
then $a<a'$ and $P(a,a')=i$.
\item[d)]$G'$ has no heavy path.
\end{description}
\end{lemma}

\proof The definition of $E'$ immediately gives a).

For b) note that all eligible edges incident to $y\in B$ have
$a_y$ for class.

For c) let $x\in A$ be the common vertex of the two edges and apply Lemma
\ref1 for $a_x$, $a$, and $a'$.

Finally d) follows since if a walk of $G'$ starts in $B$ then its
coloring $(c_1,c_2,c_3)$ must satisfy that $c_1$ and $c_2$ have
different class by c) but $c_2$ and $c_3$ have the same class
by b), so $c_2<c_1\le c_3$ is impossible. \qed

\begin{lemma}\label3 Let $G=(A,B,E)$ be a bipartite graph with a
proper edge coloring. Lexicographic
thinning produces a subgraph $G'$ with no slow walk.
\end{lemma}

\proof Suppose
$v_0v_1v_2v_3v_4$ is a walk in $G'$ starting at $v_0\in B$ and let its
coloring be $(c_1,c_2,c_3,c_4)$. Assume $c_1>c_2<c_3<c_4$. We show that
$c_1>c_4$, so this walk is not slow. By Lemma \ref2/a,b, as $c_2$ and $c_3$
are colors of edges incident to $v_2\in B$, they have the same class, but
they have different types: $c_2=(a,i,z)$, $c_3=(a,j,s)$ with $i\ne
j$. We use lexicographic ordering, so $c_2<c_3$ implies $i<j$. Both
$c_1$ and $c_2$ are colors of edges incident to $v_1\in A$, so by
Lemma \ref2/c their classes are different. Since $c_1>c_2$ we have
$b>a$ for the class $b$ of $c_1$. Still by Lemma \ref2/c $P(a,b)=i$. Similarly,
$c_3$ and $c_4$ are colors of distinct edges in $E'$ incident to $v_3\in A$,
so they have different classes. As $c_3<c_4$ we have $c>a$ for the
class $c$ of $c_4$. We have $P(a,c)=j$. By Lemma \ref1 we have $b>c$. This
proves $c_1>c_4$ as claimed. \qed

\begin{lemma}\label4
Let $G=(A,B,E)$ be a bipartite graph with a
proper edge coloring. Reversed thinning produces a subgraph
$G'$ with no fast walk.
\end{lemma}

\proof Suppose
$v_0v_1v_2v_3v_4$ is a walk in $G'$ with coloring
$(c_1,c_2,c_3,c_4)$. Assume $c_1>c_2<c_3<c_4$. We show that $c_1<c_4$,
so this walk is not fast. First assume the walk starts at $v_0\in A$. As
$G'$ does not contain a heavy path, $v_3v_2v_1v_0$ is not
heavy, so $c_1<c_3$. This implies $c_1<c_4$ as claimed.

Now assume $v_0\in
B$. Just as in the proof of the previous lemma, $c_2$ and $c_3$ are colors of
edges incident to $v_2\in B$, so they have the same class, but they have
different types: $c_2=(a,i,z)$, $c_3=(a,j,s)$ with $i\ne j$. We
use the reversed ordering, so $c_2<c_3$ implies $i>j$. Both $c_1$ and
$c_2$ are colors of edges incident to $v_1\in A$, so their classes are
different. Since $c_1>c_2$ we have $b>a$ for the class $b$ of
$c_1$ and $P(a,b)=i$. Similarly, $c_3$ and $c_4$ are colors of distinct edges
in $E'$ incident to $v_3\in A$, so they have different classes. As
$c_3<c_4$ we have $c>a$ for the class $c$ of $c_4$ and $P(a,c)=j$. By Lemma
\ref1 we have $c>b$. This proves $c_1<c_4$ as claimed. \qed

Below we estimate the number of edges in $E'$. Recall that both thinning
procedures are randomized. We can show that the subgraphs they produce
have a large expected number of edges. We did not make any effort to
optimize for the constant in this lemma.

\begin{lemma}\label5 Let $G=(A,B,E)$ be a bipartite graph with a
$d$ edge coloring. Let $G'=(A,B,E')$ be the result of either the
lexicographic or the reversed thinning. We have
$$\E[|E'|]\ge{t-1\over240d}|E|\ge{\log d\over480d}|E|.$$
\end{lemma}

\proof We compute the
probability for a fixed edge $e=(x,y)\in E$ to end up in $E'$. For this we
break down the random process producing $E'$ into three phases. In the
first phase we select $F$. With $F$ the color $\chi(e)$ is identified with an
element of $H$, most importantly, the type of $e$ is fixed. In the second
phase we select $a_x$ and $a_y$ uniformly at
random. These choices determine if $e$ is eligible. If $e$ is not
eligible then $e\notin E'$. So in the third phase we consider $F$,
$a_x$, and $a_y$ fixed and assume $e$ is eligible. We select the
random values $a_z$ for vertices $z\ne x,y$. This effects if other
edges are eligible and if $e\in E'$. Here is the detailed calculation:

Let $e\in E$ have the color $\chi(e)=k\in\{1,2,\ldots, d\}$. The choice
of $F$ in the first phase determines $F(k)=(a,i,z)\in H$. By the construction
of $F$, if we call $a'$ the last $t-1$ bits of $a$ then $(a',i,z)$ is
uniformly distributed among all its possible values. In particular,
the probability that $e$ becomes a type $i$ edge is exactly
$$\P[e\hbox{ is of type }i]=\frac{2^{t+i-2}}{2^{2t-1}-2^t}
=\frac{2^i}{2^{t+1}-4}.$$ 

For phase two we consider the function $F$ identifying colors with elements
of $H$ fixed. Consider an edge $e=(x,y)$ of
color $(a,i,z)\in H$. This edge is eligible if $a_y=a<a_x$ and
$P(a,a_x)=i$. This determines $a_y$ and the first $i$ bits of $a_x$. Recall
that by the definition of $H$ the string $a$ has $0$ in position $i$. Thus,
the probability that the edge $e$ of type $i$ is eligible is exactly
$2^{-t-i}$.

Assume for the third phase that $e$
is eligible. Consider another edge $e'=(x,y')\in E$ with color
$\chi(e')=(a',i,z')$ of type $i$. If $a$ and $a'$ do not agree in the first
$i$ positions, then $e'$ is not eligible. If they agree in the first $i$
positions, then $e'$ is eligible if and only if $a_{y'}=a'$, so with
probability $2^{-t}$. Let $k'_e$ be the number of edges $(x,y')$ of type $i$
with the first $i$ digits of their class agreeing with $a$, but not counting
$e$ itself. We have $k'_e<2^t$.

Consider now an edge $e''=(x'',y)\in E$ with color
$(a'',i,z'')$ such that $e''\ne e$. As $e$ is eligible, $e''$ can only
be eligible if $a''=a$. If $a''=a$ then $e''$ is eligible if and only
if $a_x$ and $a_{x''}$ agree in the first $i$ digits. This happens with
probability $2^{-i}$. For the number $k''_e$ of the edges $(x'',y)\ne e$ of
type $i$ and class $a$ we have $k''_e<2^i$.

In phase three the
eligibility of all these edges $e'$ and $e''$ adjacent with $e$ are
independent events.

Still consider the function $F$ fixed. The total probability for an edge $e$
of type $i$ to be in $E'$ is
$$\begin{array}{rcl}\P[e\in E'|F]
&=&2^{-t-i}(1-2^{-t})^{k'_e}(1-2^{-i})^{k''_e}\\
&\ge&2^{-t-i}(1-2^{-t})^{2^t-1}(1-2^{-i})^{2^i-1}\\
&>&2^{-t-i}/7.5.\end{array}$$

The total probability of $e\in E'$ can be calculated from the
distribution of its type and the above conditional probability
depending on its type:
$$\P[e\in E']>\sum_{i=2}^t\frac{2^i}{2^{t+1}-4}\cdot\frac{2^{-t-i}}{7.5}
>\frac{t-1}{15\cdot2^{2t}}>\frac{t-1}{240d}.$$
The expected number of edges in $E'$ is then
$$\E[|E'|]>\frac{t-1}{240d}|E|\ge\frac{\log d}{480d}|E|.$$
\qed

\begin{theorem}\label6
Let $G=(A,B,E)$ be a bipartite graph
with a proper edge coloring. There exists a subgraph
$G'=(A,B,E')$ of $G$ without a slow walk and with $|E'|>\frac{\log
d}{480d}\cdot|E|$. Similarly, there exists a subgraph
$G''=(A,B,E'')$ of $G$ without a fast walk and with $|E''|>\frac{\log
d}{480d}\cdot|E|$.
\end{theorem}

\proof By Lemmas \ref3 and \ref4 the results of the lexicographic and reversed
thinnings avoid the slow and fast walks, respectively. There exists an
instance of the random choices with the size of $E'$ being at least
its expectation given in Lemma \ref5. This proves the theorem. \qed

\subsection{Longer forbidden walks}

Here we generalize the concept of fast and slow walks to longer
walks. Consider a bipartite graph $G=(A,B,E)$ with a proper edge
coloring. For $k\ge2$ we call a walk of length $2k$ in $G$ a {\em
$k$-fast walk} if its coloring $(c_1,\ldots,c_{2k})$ satisfies
$c_1>c_2>\ldots>c_k<c_{k+1}<c_{k+2}<\ldots<c_{2k}$ and $c_1\ge
c_{2k}$. For $k\ge2$ we call a walk of length $2k$ in $G$ a {\em $k$-slow
walk} if its coloring $(c_1,\ldots,c_{2k})$ satisfies the following:
$c_{2j-1}>c_{2j}$ for
$1\le j\le k/2$; $c_{2j}<c_{2j+1}$ for $1\le j<k/2$; $c_{2j-1}<c_{2j}$
for $k/2<j\le k$; $c_{2j}>c_{2j+1}$ for $k/2\le j<k$; and finally
$c_1\ge c_{2k}$. If a $k$-slow walk starts in the vertex set $B$ we call it a
$(k,B)$-slow walk, otherwise it is a $(k,A)$-slow walk.

Notice that $2$-fast walks are the fast walks
and $(2,B)$-slow walks are the slow walks with their orientation
reversed. For the coloring $(c_1,\ldots,c_{2k})$ of a
$k$-fast walk $c_j$ is in between $c_{j-1}$ and $c_{j+1}$ for all
$1<j<2k$, $j\ne k$, while $c_k$ is the smallest color on this
list. For the coloring $(c_1,\ldots,c_{2k})$ of a $k$-slow walk the
situation is reversed: the only index $1<j<2k$ with $c_j$ being in
between $c_{j-1}$ and $c_{j+1}$ is the index $j=k$.

In order to apply the lexicographic and reversed thinning recursively
we have to change the coloring of the subgraph. Let $G=(A,B,E)$ be a
bipartite graph with a proper edge coloring given by
$\chi:E\to\{1,\ldots,d\}$. Let $G'=(A,B,E')$ the result of the
lexicographic or the reversed thinning of $G$. Recall that the edges
in $E'$ have a type $2\le i\le t$ with $t=\lceil(\log
d)/2\rceil+1$. The {\em type edge coloring} of $G'$ is the map
$\chi':E'\to\{1,\ldots,t-1\}$ defined by $\chi'(e)=t+1-i$ for an edge $e\in
E'$ of type $i$. By Lemma \ref2/a $\chi'$ is a proper edge coloring of $G'$.

\begin{lemma}\label7
Let $G=(A,B,E)$ be a bipartite graph with proper
edge coloring given by $\chi:E\to\{1,\ldots,d\}$. Let $G'=(A,B,E')$ the
result of the lexicographic thinning of $G$. Let $\chi'$ be the type edge
coloring of $G'$ and let $k\ge2$. If a subgraph $G''=(A,B,E'')$ of $G'$
with its edge coloring given by $\chi'$ has no $(k',A)$-slow walk for
$2\le k'<k$ then $G''$ with its edge coloring given by $\chi$ has no
$(k,B)$-slow walk.
\end{lemma}

Notice that the $k=2$ case of this lemma gives a second proof of
Lemma~\ref3.

\proof Let $W=v_0v_1\ldots v_{2k}$ be a walk in $G''$ starting at
$v_0\in B$ and let its $\chi$-coloring be $(c_1,c_2,\ldots,c_{2k})$. Assume
that $c_i>c_{i+1}$ or $c_i<c_{i+1}$ for $1\le i<2k$ as required in the
definition of a $(k,B)$-slow walk. We need to show $c_1<c_{2k}$.

We identify the colors of $\chi$ with the triplets $(a,i,z)\in H$ as in
the definition of lexicographic thinning. We let
$c_j=(a_j,i_j,z_j)$. The $\chi'$-coloring of $W$ is
$(t+1-i_1,\ldots,t+1-i_{2k})$. We have $i_j\ne i_{j+1}$ for $1\le j<2k$.

For $1\le j<k$ the colors $c_{2j}$ and $c_{2j+1}$ are colors
of distinct edges incident to $v_{2j}\in B$, so by Lemma \ref2/b their
class is the same: $a_{2j}=a_{2j+1}$. We consider lexicographic thinning, so
the order between $c_{2j}$ and $c_{2j+1}$ is the same as the order between
their types: $i_{2j}$ and $i_{2j+1}$.
For $1\le j<k/2$ we have $i_{2j}<i_{2j+1}$ but for
$k/2\le j<k$ we have $i_{2j}>i_{2j+1}$.

For $1\le j\le k$ the colors
$c_{2j-1}$ and $c_{2j}$ are colors of edges incident to $v_{2j-1}\in
A$. By Lemma \ref2/c the classes of these colors do not agree, and the
ordering between the classes, between the types, and between the
colors themselves are the same. Thus, for $1\le j\le k/2$ we have
$a_{2j-1}>a_{2j}$ and $i_{2j-1}>i_{2j}$. For $k/2<j\le k$ we have
$a_{2j-1}<a_{2j}$ and $i_{2j-1}<i_{2j}$. Also by Lemma \ref2/c for all
$1\le j\le k$ we have $P(a_{2j-1},a_{2j})=\min(i_{2j-1},i_{2j})$.

The sequence $a_1,a_2,\ldots,a_k$ is monotone
decreasing and it changes only in every other step. The first
positions of change between distinct consecutive elements are
$i_2,i_4,\ldots,i_{2\lfloor k/2\rfloor}$. So we have $a_1>a_k$ and
$P(a_1,a_k)=\min(S_1)$ for the set $S_1=\{i_2,i_4,\ldots,i_{2\lfloor
k/2\rfloor}\}$.

Similarly, $a_k,a_{k+1},\ldots,a_{2k}$ is monotone
increasing and it changes only in every other step. The first
positions of change between distinct consecutive elements are $i_{2\lfloor
k/2\rfloor+1},\ldots,i_{2k-3},i_{2k-1}$. So we have $a_k<a_{2k}$ and
$P(a_k,a_{2k})=\min(S_2)$ for the set $S_2=\{i_{2\lfloor
k/2\rfloor+1},\ldots,i_{2k-3},i_{2k-1}\}$.

Let us consider an arbitrary value $1\le j<k/2$ and let
$2\le k'=k-2j+1<k$. Consider the $2k'$ long middle portion $W'$ of the
walk $W$: let $W'=v_{2j-1}v_{2j}\ldots v_{2k-2j+1}$. This is a walk of
length $2k'$ in $G''$ starting at $v_{2j-1}\in A$. By our assumption on $G''$
this is not a $(k',A)$-slow walk if considered with the
coloring $\chi'$. But the $\chi'$-coloring
of $W'$ is $(t+1-i_{2j},t+1-i_{2j+1},\ldots,t+1-i_{2k-2j+1})$ and the
consecutive values in this list compare as required for a $(k',A)$-slow
walk. Therefore, we must have $t+1-i_{2j}<t+1-i_{2k-2j+1}$.

We have just proved $i_{2j}>i_{2k-2j+1}$ for $1\le j<k/2$. For even $k$ and
$j=k/2$ the same formula compares the types of two consecutive edges of $W$
and we have already seen its validity in that case too. For every element of
the set $S_1$ we have just found a smaller element of the set
$S_2$. Therefore, $\min(S_1)>\min(S_2)$. Using that $a_{2k}>a_k$ and
$P(a_1,a_k)=\min(S_1)>\min(S_2)=P(a_k,a_{2k})$ Lemma \ref1 gives
$a_1<a_{2k}$. This implies $c_1<c_{2k}$ and finishes the proof of the
lemma. \qed

\begin{lemma}\label8
Let $G=(A,B,E)$ be a bipartite graph with proper
edge coloring given by $\chi:E\to\{1,\ldots,d\}$. Let $G'=(A,B,E')$ the
result of the reversed thinning of $G$. Let $\chi'$ be the type edge
coloring of $G'$ and let $k\ge2$. If a subgraph $G''=(A,B,E'')$ of $G'$
with its edge coloring given by $\chi'$ has no $(k',A)$-slow walk for
$2\le k'<k$ then $G''$ with its edge coloring given by $\chi$ has no
$k$-fast walk.
\end{lemma}

Notice that the $k=2$ case of this lemma gives a second proof of
Lemma~\ref4.

\proof The proof of this lemma is very similar to that of Lemma \ref7.

Let $W=v_0v_1\ldots v_{2k}$ be a walk in $G''$ and let its
$\chi$-coloring be $(c_1,c_2,\ldots,c_{2k})$. Assume that
$c_1>c_2>\ldots>c_k<c_{k+1}<c_{k+2}<\ldots<c_{2k}$ as required in the
definition of a $k$-fast walk. We need to show $c_1<c_{2k}$. Instead, we prove
the slightly stronger statement that the class of $c_1$ is smaller than
the class of $c_{2k}$. We first do that for walks starting in $B$:
assume that $v_0\in B$.

We identify the colors of $\chi$ with the triplets $(a,i,z)\in H$ as in
the definition of reverse thinning. We let
$c_j=(a_j,i_j,z_j)$. Note that the $\chi'$-coloring of $W$ is
$(t+1-i_1,\ldots,t+1-i_{2k})$. We have $i_j\ne i_{j+1}$ for $1\le j<2k$.

For $1\le j<k$ the colors $c_{2j}$ and $c_{2j+1}$ are colors
of distinct edges incident to $v_{2j}\in B$, so by Lemma \ref2/b their
classes are the same: $a_{2j}=a_{2j+1}$. Thus the order between $c_{2j}$
and $c_{2j+1}$ is determined by the order between $i_{2j}$ and
$i_{2j+1}$, but as we use the reversed ordering in $H$ the order
between $c_{2j}$ and $c_{2j+1}$ is reversed compared to the order
between $i_{2j}$ and $i_{2j+1}$. Specifically, for $1\le j<k/2$ we
have $i_{2j}<i_{2j+1}$ and for $k/2\le j<k$ we have
$i_{2j}>i_{2j+1}$. For $1\le j\le k$ the colors $c_{2j-1}$ and
$c_{2j}$ are colors of edges incident to $v_{2j-1}\in A$. By Lemma \ref2/c
the classes of these colors do not agree, and the ordering between the
classes, between the types, and between the colors themselves are the
same. Thus, for $1\le j\le k/2$ we have $a_{2j-1}>a_{2j}$ and
$i_{2j-1}>i_{2j}$. For $k/2<j\le k$ we have $a_{2j-1}<a_{2j}$ and
$i_{2j-1}<i_{2j}$. Also by Lemma \ref2/c for $1\le j\le k$ we have
$P(a_{2j-1},a_{2j})=\min(i_{2j-1},i_{2j})$.

At this point we have the same ordering of the classes and types of
the coloring of $W$ as in the proof of Lemma \ref7. We also have the same
assumption that $G''$ with the edge coloring $\chi'$ has no $(k',A)$-slow
walk for $2\le k'<k$. So we arrive to the same conclusion $a_1<a_{2k}$
with an identical proof.

We finish the proof of the lemma by considering the alternative case
when $W$ starts in $A$. Now $c_1$ and $c_2$ are colors of edges
sharing a vertex $v_1\in B$, so by Lemma \ref2/b their classes are
equal. Similarly, the classes of $c_{2k-1}$ and $c_{2k}$ are equal, so
it is enough to prove that the class of $c_2$ is smaller than the
class of $c_{2k-1}$. For $k=2$ this follows directly from Lemma
\ref2/c. For $k>2$ the walk $W'=v_1v_2\ldots v_{2k-1}$ is exactly the type
of walk we considered for $k_0=k-1$. As it starts in $B$ we have
already proved that the class of its first edge is smaller than the
class of its last edge. This finishes the proof of the case of a walk
starting in $A$ and also the proof of Lemma \ref8. \qed

Lemmas \ref7 and \ref8 set the stage to use the thinning procedures
recursively to get subgraphs avoiding $(k,B)$-slow or $k$-fast walks.
In a single application of either thinning procedure the number $d$ of
colors in the original coloring is replaced by $t-1=\lceil\log
d/2\rceil$ colors in the type coloring. Here $4(t-1)>\log(4d)$, so
after $k$ recursive calls we still have more than $\log^{(k)}(4d)/4$
colors, where $\log^{(k)}$ stands for the $k$ times iterated log
functions. (Of course, this only makes sense if $\log^{(k)}(4d)>2$. Otherwise
we can get stuck, as neither thinning procedure is defined in the pathetic
case of $d=1$ colors.) Making optimal random choices we may assume that each
thinning procedure yields at least the expected number of
edges. Thus, the ratio of the number of edges and the number of colors
decreases by at most a factor of $240$ in each iteration. Clearly,
the only interesting case is when the original average degree was
$\Theta(d)$ in which case the average degree after $k$ iterations
remains $\Theta(\log^{(k)}d)$. The constant of proportionality depends
on $k$.

\begin{theorem}\label9
Let $G=(A,B,E)$ be a bipartite graph
with a $d$ edge coloring and let $k\ge2$. There exists a subgraph
$G'=(A,B,E')$ of $G$ without a $(k',B)$-slow walk for any $2\le k'\le k$
and with $|E'|>{\log^{(k-1)}d\over4\cdot240^{k-1}d}|E|$. Similarly, there
exists a subgraph $G''=(A,B,E'')$ of $G$ without a $k'$-fast walk for
any $2\le k'\le k$ and with
$|E''|>{\log^{(k-1)}d\over4\cdot240^{k-1}d}|E|$.
\end{theorem}

\proof We apply the thinning procedures recursively. First we use
lexicographic and reversed thinning to obtain subgraphs $G_1$ and $G_2$ of
$G$, respectively. We make sure these graphs have at least as many edges as
the expected number given in Lemma \ref5. If $k=2$ we are done, $G'=G_1$
and $G''=G_2$ satisfy the conditions of the theorem. Otherwise we consider
$G_1$ and $G_2$ with
the type edge coloring. We find recursively their subgraphs $G'$ and $G''$,
respectively, avoiding $(k',A)$-slow walks for $2\le k'\le k-1$. This can be
done because the sides $A$ and $B$ play symmetric roles. Finally, we apply
Lemmas \ref7 and Lemma \ref8 to see that the subgraphs $G'$ and $G''$, if
considered with the original edge coloring of $G$, avoid all walks required in
the theorem. The number of edges guaranteed in the subgraphs is calculated in
the paragraph preceding the theorem and is at least the stated bound. \qed

\subsection{$k$-flat graphs}

In this subsection we establish that removing a linear number of edges
from a $k$-fast walk free graph the resulting graph has special
structural properties. We note here that the recursive thinning
construction that we used to arrive at $k$-fast walk free graphs results
in a graph that itself is $k$-flat as defined below. We chose however
to keep the inductive part of the proof simple and concentrated only
on $(k,B)$-slow and $k$-fast walks. We derive the more complicated
properties from these simpler ones. Note that in this subsection we
do not use that our graphs are bipartite.

Let $G$ be a graph and $\chi$ a proper edge coloring of $G$. We define
the {\em shaving} of the graph $G$ to be the subgraph obtained from
$G$ by deleting the edge with the largesr color incident to every
(non-isolated) vertex. Clearly, we delete at most $n$ edges, where $n$
is the number of vertices in $G$. We define the {\em $k$-shaving} of
$G$ to be the subgraph obtained from $G$ by repeating the shaving operation
$k$ times. Clearly, we delete at most $kn$ edges for a $k$-shaving.

Let $W$ be a walk of length $m$ in a properly edge colored graph $G$,
and assume its coloring is $(c_1,\ldots,c_m)$. We define the {\em height
function} $h_W$ from $\{1,\ldots,m\}$ to the integers recursively
letting $h_W(1)=0$ and
$$h_W(i+1)=\left\{\begin{array}{lrl}h_W(i)+1&\hbox{~~~if }&c_{i+1}>c_i\\
h_W(i)-1&\hbox{if }&c_{i+1}<c_i\end{array}\right.$$
for $1\le i<m$. Note $h_W(i)+i$ is always odd. This function considers
how the colors of the edges in the walk change, in particular, how
many times the next color is larger and how many times it is smaller
than the previous color.

We call a graph $G$ with proper edge coloring {\em $k$-flat} if the
following is true for every walk $W$ in $G$. Let $m\ge2$ be the length of
$W$, let $(c_1,\ldots,c_m)$ be the coloring of $W$ and assume that the
height function satisfies $h_W(i)<0$ for $2\le i\le m$. If $m\le 2k+1$
or $h_W(i)\ge-k$ for all $i$ then $c_1>c_m$.

\begin{lemma}\label{10}
Let $G$ be properly edge colored graph. Let
$k\ge1$ and assume $G$ has no $k'$-fast walk for $2\le k'\le k$. Then
the $(k-1)$-shaving $G'$ of $G$ is $k$-flat.
\end{lemma}

\proof We prove the following slightly stronger statement by induction on
$m$. Let $W=v_0\ldots v_m$ be a walk of length $m$ in $G$ with
coloring $(c_1,\ldots,c_m)$. Let $1\le j\le m$ be the largest index such
that $h_W(j)=1-j$. Assume the walk $v_jv_{j+1}\ldots v_m$ is in the
$(k-1)$-shaving $G'$ of $G$. Also assume that $h_W(i)<0$ for $2\le
i\le m$. If $m\le 2k+1$ or $h_W(i)\ge-k$ for all $i$ then we claim
$c_1>c_m$. This statement is stronger than Lemma \ref{10} since it allows for
the initial decreasing segment of $W$ be outside $G'$.

If $j=m$ the statement of the claim is obvious from the definition of
the height function. This covers the $m=2$ and $m=3$ base cases. Let
$m\ge4$ and assume the statement is true for walks of length $m-1$
and $m-2$.

If $j>k+1$ we have $h_W(j)=1-j<-k$ so we must have
$m\le2k+1$. Consider the walk $W'=v_1\ldots v_m$ of length $m-1$. We
have $h_{W'}(i)=1-i<0$ for $2\le i<j$ and $h_{W'}(i)\le
h_{W'}(j-1)+(i-(j-1))=3+i-2j<0$ for $j\le i\le m-1$. Thus the
inductive hypothesis is applicable and we get $c_1>c_2>c_m$.

Finally consider the $j\le k+1$ case. As the trivial $j=m$ case was
already dealt with we also assume $j<m$. Clearly, $h_W(2)<0$ implies
$j\ge2$. We chose a $w_0\ldots w_{j-2}$
walk in $G$ ending at $w_{j-2}=v_j$ and with coloring
$(c_1',\ldots,c_{j-2}')$ satisfying
$c_1'>c_2'>\ldots>c_{j-2}'>c_{j+1}$. This is possible since the edge
$v_jv_{j+1}$ is in the $(k-1)$-shaving $G'$ of $G$, so we can find the
edge $w_{j-3}w_{j-2}$ in the $(k-2)$-shaving of $G$, $w_{j-4}w_{j-3}$
in the $(k-3)$-shaving, and so on. We must have $c_1>c_1'$ as otherwise
$w_0w_1\ldots w_{j-3}v_jv_{j-1}\ldots v_0$ is a $(j-1)$-fast walk
and no such walk exists in $G$. Now consider the walk $W'=w_0w_1\ldots
w_{j-3}v_jv_{j+1}\ldots v_m$. This is a walk of length $m-2$ and
satisfies $h_{W'}(i)=1-i$ for $1\le i\le j-1$ and $h_{W'}(i)=h_W(i+2)$ for
$j-1\le i\le m-2$. All requirements of the inductive hypothesis are
satisfied, so we have $c_1'>c_m$. Thus $c_1>c_1'>c_m$ as claimed. \qed

\begin{corollary}\label{corol}
Let $G=(A,B,E)$ be a bipartite graph
with a $d$ edge coloring and let $k\ge2$. There exists a $k$-flat
subgraph $G'=(A,B,E')$ of $G$ with
$|E'|>{\log^{(k-1)}d\over4\cdot240^{k-1}d}|E|-(k-1)(|A|+|B|)$.
\end{corollary}

\proof Combine Theorem~\ref9 nwith Lemma~\ref{10} and the fact that
$(k-1)$-shaving keeps all but at at most $(k-1)(|A|+|B|)$ edges of $G$. \qed
\smallskip

The final lemma in this section is a simple but useful observation on $k$-flat
graphs. It can also be stated for longer walks with height function
bounded from below, but for simplicity we restrict attention to short
walks.

\begin{lemma}\label{11}
Let $k\ge1$ and let $G$ be a properly edge
colored $k$-flat graph. Let $W=v_0\ldots v_m$ be a walk in $G$ of
length $m\le2k+1$ with coloring $(c_1,\ldots,c_m)$. If $c_1\ge c_i$
for all $1\le i\le m$ then $h_W(i)\le0$ for all $1\le i\le m$.
\end{lemma}

\proof We prove the contrapositive statement. Assume $h_W(i)>0$ for
some $1\le i\le m$ and let $i_0$ be the smallest such index. Clearly,
$i_0\ge2$, $h_W(i_0)=1$ and for the walk $W'=v_{i_0}v_{i_0-1}\ldots
v_0$ we have $h_{W'}(i)=h_W(i_0-i+1)-1<0$ for $1<i\le i_0$. So by the
definition of $k$-flatness we have $c_{i_0}>c_1$. \qed

\section{Locally plane graphs}
\label{s3}

Locally plane graphs were introduced in the paper \cite{PPTT} (though the name
appears first in this paper). That paper gives a
simple construction for $3$-locally plane graphs. We recall (a simplified
version of) the construction as it is our starting point.

\subsection{Construction of $3$-locally plane graphs in \cite{PPTT}}

Let $d\ge1$ and consider the orthogonal projection of (the edge graph of)
the $d$ dimensional hypercube into the plane. A suitable projection of the
``middle layer'' of the hypercube provides the $3$-locally plane graph. Here
is the construction in detail:

Let $d\ge1$ be fixed and set $b=\lfloor d/2\rfloor$. The bit at position $i$
in $x\in\{0,1\}^d$ (the $i$\/th coordinate) is denoted by $x_i$ for $1\le i\le
d$. We let $A=\{x\in\{0,1\}^d\mid\sum_{i=1}^dx_i=b\}$ and
$B=\{x\in\{0,1\}^d\mid\sum_{i=1}^dx_i=b+1\}$. The abstract graph
underlying the geometric graph to be constructed is $G_d=(A,B,E)$ with
$(x,y)\in E$ if $x\in A$, $y\in B$ and $x$ differs from $y$ in a
single position. This is the middle layer of the $d$ dimensional
hypercube. We define the edge coloring $\chi:E\to\{1,\ldots,d\}$ that colors
an edge $e=(x,y)\in E$ by the unique position $\chi(e)=i$ with $x_i\ne
y_i$. Notice that this is a proper edge coloring. The number of vertices is
$n=|A|+|B|={d\choose b}+{d\choose b+1}\le2^d$, the number of edges is
$|E|={d\choose b}(d-b)>nd/4$. The average degree is greater than $d/2\ge\log
n/2$.

To make the abstract graph $G_d$ into a geometric graph we project the
hypercube into the plane. We give two possible projections here. The first is
more intuitive and it is closer to the actual construction in \cite{PPTT}. We
let $a_i=(10^i,i\cdot10^i)$ for $1\le i\le d$ and use this vector as the
projection of the edges of color $i$.
We will use that among the vectors $a_i$ higher
index means higher slope and much greater length. We identify
the vertex $x\in A\cup B$ with the point $P_x=\sum_{i=1}^dx_ia_i$. The
edges are represented by the straight line segment connecting their
endpoints.

We give the second construction to obtain a graph where all edges are very
close in length and direction. Let $0<\epsilon<10^{-d}$ be arbitrary and consider
the vectors $b_i=(1+10^i\epsilon,\epsilon^{d+1-i}(1+10^i\epsilon))$ and identify
a point $x\in A\cup B$ with $Q_x=\sum_{i=1}^dx_ib_i$.

It is easy to verify that we get a geometric
graph in both cases (i.e., the vertices are mapped to distinct points and
no edge passes through a vertex that is not its endpoint). Note that
edges of color $i$ are all translates of the same vector $a_i$ or $b_i$. We do
not introduce separate notations for the two geometric graphs constructed this
way as they will only be treated separately in the proof of
Lemma~\ref{12}, where we refer to them as the first and the second realization
of $G_d$.

\subsection{Self-intersecting paths in $G_d$}

In \cite{PPTT} a graph very similar to $G_d$ was shown to be $3$-locally
plane. Here we do more, we analyze
all self-intersecting paths of $G_d$ as follows.

\begin{lemma}\label{12}
Let $W$ be a walk in $G_d$ with coloring $(c_1,\ldots,c_m)$ satisfying $c_1\ge
c_m$. Assume that $W$ and all its non-empty subwalks have a unique edge of maximal
color. The first and last edges of $W$ cross in either geometric
realization if and only if $m$ is even and
there is an odd index $1<j<m$ satisfying $c_1>c_m>c_j\ge c_i$ for all $1<i<m$.
\end{lemma}

\proof Let $W=v_0\ldots v_m$. Note that the first and the last edges cross if
and only if $v_0$ and $v_1$ are on different sides of the line $\ell$ through
$v_{m-1}$ and $v_m$ and similarly $v_{m-1}$ and $v_m$ are on different sides
of the line $\ell'$ through $v_0$ and $v_1$. To analyze such separations
consider the projection $\pi_i$ to the $y$ axis parallel to edges of color $i$.

Let us consider the first realization of $G_d$ with the
vectors $a_i$. We have $\pi_i(x,y)=y-ix$ and the projection of the vector
$a_j$ is of length $|i-j|10^j$. Thus higher colored edges map to longer
intervals (except color $i$ itself). Under the projection $\pi_{c_1}$ the
direction of the highest colored edge in the walks $v_1\ldots v_{m-1}$
(respectively, $v_1\ldots v_m$) determines which
side $v_{m-1}$ (respectively $v_m$) lies of the line $\ell'$. Indeed this
highest color cannot be $c_1$, so the projections of the other edges will be
much shorter and by the
unique maximal color property we see that the walk contains at most
$2^{k-1}$ edges having the $k$\/th largest color, so these shorter projections
cannot add up to be more than the projection of the largest edge.

We can only have $v_{m-1}$ and $v_m$ lying on different sides of $\ell'$ if
these edges of maximal color are distinct, thus we must have $c_m>c_i$ for all
$1<i<m$. From $c_1\ge c_m$ and the unique maximal color property we have
$c_1>c_m$. Taking $c_j$ to maximize $c_i$ for $1<i<m$ (this is unique again)
we have $c_1>c_m>c_j\ge c_i$ for all $1<i<m$. 

It is left to prove that the first and last edges of $W$ cross if and only if
$m$ is even and $j$ is odd.

To prove this claim one has to use that $G_d$ is
bipartite with vertex sets $A$ and $B$, and every edge of
color $c$ is a translate of the vector $a_c$ with its head in $B$ and
tail in $A$. Thus, the vector $v_{i-1}v_i$ is either $a_{c_i}$ or
$-a_{c_i}$ depending on the parity of $i$. Which sides of $\ell'$ $v_{m-1}$
and $v_m$ lie is determined by the projections of the $j$\/th and last edge,
so they are on opposite sides if $j$ and $m$ have different
parities. Similarly, the sides of $\ell$ on which $v_0$ and $v_1$ lies is
detemined by the $\pi_{c_m}$ projection of the first and $j$\/th edges, but as
we have $c_1>c_m>c_j$ they are on different sides if $1$ and $j$ has the same
parity. See Figure~1 for a rough depiction
of all four cases. This finishes the proof of the claim and the part of the lemma
regarding the first realization of $G_d$.

\begin{figure}[ht]
\begin{center}
\scalebox{0.6}{\includegraphics{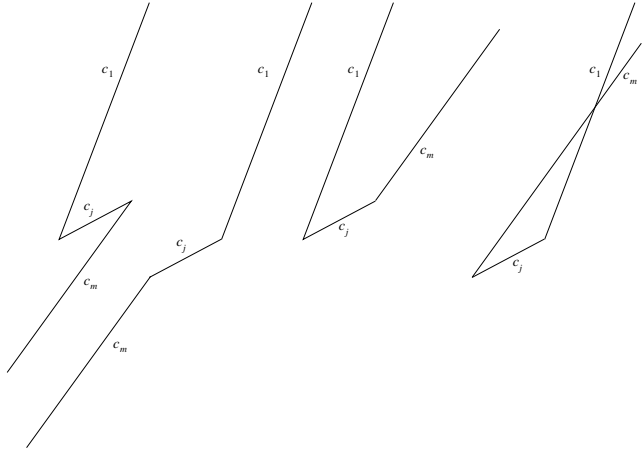}}
\smallskip
\caption{Rough picture of $W$ in the first realization of $G_d$ in the four cases according to
  the parities of $m$ and $j$. The scale is set to match the edge $v_{j-1}v_j$
  of color $c_j$, shorter edges are approximated by zero, and only initial
  segments of the two longer edges are depicted.}
\end{center}
\end{figure} 

The proof for the second realization of $G_d$ as a geometric graph (involving
the vectors $b_i$) is slightly more complicated. We have
$\pi_i(x,y)=y-\epsilon^{d+1-i}x$ and
$\pi_i(b_j)=\epsilon^{d+1-j}(1+10^j\epsilon)-10^j\epsilon^{d+2-i}-\epsilon^{d+1-i}$.
The $\epsilon^{d+1-i}$ terms alternate in sign in the projection on the edges
and cancel completely for a walk of even length. For a walk of odd length a single
such term remains but it is dominated by the other terms if the walk has an
edge with color above $i$. If, however, no such edge exists, the remaining
uncanceled $\epsilon^{d+1-i}$ term dominates the other terms in the
projection. Thus, if the projection of the
unique edge with the largest color of a walk has color $c>i$ or if $c<i$ but the walk has even
length, then the sign of the $\pi_i$ projection of the edge with the largest color determines
the sign of the projection of the entire walk. But in case $c<i$ and the walk has
odd length, then the non-canceling $\epsilon^{d+1-i}$ term determines the sign.

Let the edge $v_{j-1}v_j$ be the one with the unique largest color in the walk
$v_1\ldots v_{m-1}$. A case analysis of the parities of $j$, $m$ and whether
$c_j>c_1$ or $c_j>c_m$ hold show that the first and last edges of
$W$ cross if and only if $c_m>c_j$, $j$ is odd and $m$ is even -- as claimed
in the lemma. \qed

We call $v_m\ldots v_0$ the {\em reverse} of the walk $W=v_0\ldots v_m$.

\begin{lemma}\label{13}
Let $k\ge1$ and let $G'$ be a $k$-flat subgraph of $G_d$. If the length of a
walk $W$ in $G'$ does not exceed $2k+1$, then $W$ has unique edge of largest
color.
\end{lemma}

\proof Let $W$ be a walk of length $3\le m\le2k+1$ in $G'$ with
coloring $(c_1,\ldots,c_m)$ and assume the largest color is not
unique. We may assume $c_1=c_m>c_i$ holds for all $1<i<m$, otherwise
one can take a suitable subwalk of $W$. By Lemma~\ref{11} we have
$h_W(m)\le0$. Consider $W'$ the reverse of the walk $W$. Clearly,
$h_{W'}(m)=-h_W(m)$, but also by Lemma~\ref{11} we have $h_{W'}(m)\le0$. So
we must have $h_W(m)=0$ and $m$ must be odd.

The contradiction that proves
the statement of the lemma comes from the simple observation that between two
consecutive appearances of a color in any walk of $G_d$ there always
are an even number of edges, so $m$ must be even. To see this recall that
$G_d$ is bipartite with sides $A$ and $B$ where $A$ consists of the 0-1
sequences of length $d$ with $\lfloor d/2\rfloor$ ones, while the 0-1
sequences in $B$ contain one more ones. The $A$ end of an edge of color $c$
has $0$ at position $c$, while the $B$ end of this edge has $1$ there. Along
edges of other colors bit $c$ does not change. Thus if a walk
traverses an edge of color $c$ from $A$ to $B$ say, then along the walk
bit $c$ remains $1$ until the next time the walk traverses an edge of color
$c$ and this has to be from $B$ to $A$. \qed

We note that Lemma~\ref{13} immediately implies that the girth of a $k$-flat
subgraph of $G_d$ is at least $2k+2$. This estimate can be improved by
observing that any cycle in $G_d$ has an even
number of edges of any color (one has to flip a bit even times to get
beck to the original state). In particular any cycle has at least
two occurrences of the largest color. These edges break the cycle into
two paths sharing two edges. Both path has to be of length at
least $2k+2$, the length of the cycle is at least $4k+2$.

As every properly edge colored graph is $1$-flat the $k=1$ case of the
next lemma establishes that $G_d$ is $3$-locally plane.

\begin{lemma}\label{14}
For any $k\ge1$ any $k$-flat subgraph of $G_d$ is $(2k+1)$-locally plane in either realization.
\end{lemma}

\proof Let $G'$ be a $k$-flat subgraph of $G_d$. We need to show that
no walk (or path) $W$ of length $m\le2k+1$ is self-intersecting. It is
clearly enough to show that the first and the last edges of $W$ do not
cross and we may assume that the color of the first edge is not smaller than
that of the last edge (otherwise simply consider the same walk reversed).
By Lemma~\ref{13} $W$ and all its subwalks have a unique edge
of maximal color, so Lemma \ref{12} applies. It is enough to show that the
coloring $(c_1,\ldots,c_m)$ of $W$ does not satisfy the conditions of
Lemma \ref{12}. Assume the contrary. So $m$ is even, and there is an odd
index $j$ such that $c_1>c_m>c_j\ge c_i$ for all $1<i<m$.

Consider the walk $W_1=v_mv_{m-1}\ldots v_{j-1}$ of length
$m-j+1$. Its coloring is $(c_m,\ldots,c_j)$ and $c_m$ is its largest
color. By Lemma \ref{11} $h_{W_1}(m-j+1)\le0$. In fact, as $m-j+1$ is even
$h_{W_1}(m-j+1)$ is odd, so we have $h_{W_1}(m-j+1)\le-1$.

Consider the walk $W_2=v_jv_{j-1}\ldots v_1$ of length $j-1$ and its
coloring $(c_j,\ldots,c_2)$. The largest color in $W_2$ is $c_j$, so
we have $h_{W_2}(j-1)\le0$ by Lemma \ref{11}. And again by parity
considerations $h_{W_2}(j-1)\le-1$.

It is easy to see that $h_W(m)=-1-h_{W_2}(j-1)-h_{W_1}(m-j+1)$. So we
have $h_W(m)\ge-1+1+1=1$ contradicting Lemma \ref{11}. The contradiction
proves Lemma \ref{14}. \qed

\begin{theorem}\label{15}
For any fixed $k>0$ and large enough $n$ there exists a $(2k+1)$-locally
plane graphs on $n$ vertices having at least $({\log^{(k)}n\over240^k}-k)n$
edges. Given two arbitrary disks in the plane, one can further assume that all
edges of these graph connect a vertex inside one disk to one inside the
other disk.
\end{theorem}

\proof Simply combine the results of Corollary~\ref{corol} and
Lemma~\ref{14}. If $n$ is not the size of the vertex set of $G_d$ for any $d$,
add isolated vertices to the largest $G_d$ with fewer than $n$ vertices. Use
the second realization of $G_d$ with a small enough
$\epsilon>0$ to obtain a geometric graph with all edges connecting two small
disks and apply a homothety and a rotation to get to the desired disks.
\qed

\section{Discussion on optimality of thinning}
\label{s4}

The maximum number of edges of a $3$-locally plane graph on $n$ vertices is
$\Theta(n\log n)$ as proved in \cite{PPTT}. The lower
bound is reproduced here by the $k=1$ case of Theorem \ref{15}, which is
therefore tight. The upper bound of \cite{PPTT} extends to {\em $x$-monotone
topological graphs}, i.e., when the edges are represented by curves
with the property that every line parallel to the $y$ axis intersects
an edge at most once. Without this artificial assumption on
$x$-monotonicity only much weaker upper bounds are known. For
higher values of $k$ we do not have tight results even if the
edges are straight line segments as considered in this paper. While
the number of edges in a $k$-locally plane graph constructed here deteriorates
very rapidly with the increase of $k$, the upper bound
hardly changes. In fact the only known upper bound better than
$O(n\log n)$ is for $5$-locally plane graphs: they have $O(n\log
n/\log\log n)$ edges as shown in \cite{PPTT}. Slightly better bounds are
known for geometric (or $x$-monotone topological) graphs with the
additional condition that a vertical line intersects every edge. If such a graph
is $(2k)$-locally plane for $k\ge2$, then it has
$O(n\log^{1/k}n)$ edges. The first realization of $G_d$
does not satisfy this condition, but the second one does.
Still, the lower and upper bounds
for this restricted problem are far apart: for $4$-locally plane
graphs with a cutting line the upper bound on the number of edges is
$O(n\sqrt{\log n})$ while the construction gives $\Omega(n\log\log
n)$.

Although we cannot establish that the locally plane graphs
constructed are optimal we can prove that the thinning procedure we use
is optimal within a constant factor. It follows that any $4$-locally
plane subgraph of either realization of $G_d$ has $O(n\log\log n)$ edges. This optimality
result below refers to a single step of the thinning procedure. It would be
interesting to establish a strong upper bound on the number of edges of a
$k$-flat graph for $k\ge3$.

Let us mention here that the thinning procedures described in this
paper found application in the extremal theory of $0$-$1$ matrices, see
\cite{T}, and there the result is shown to be optimal within a constant
factor. Consider an $n$ by $n$ \ 0-1 matrix that has no $2$ by $3$ submatrix
of either of the following two forms:
$$\left(\begin{array}{ccc}1&1&*\\1&*&1\end{array}\right),\hskip1cm
\left(\begin{array}{ccc}1&*&1\\{*}&1&1\end{array}\right),$$
where the $*$ can represent any entry. The maximal number of $1$ entries in
such a matrix is $\Theta(n\log\log n)$ as proved in \cite{T}. The construction
proving the lower bound is based on lexicographic thinning.

The following lemma shows that the number of edges in the subgraphs
claimed in Lemma \ref0 and Theorem \ref6 are optimal in a very strong sense:
no properly edge colored graph with significantly more edges than the
ones guaranteed by the above results can avoid heavy paths (slow or fast
walks, respectively).

\begin{lemma}\label{16}
Let $G=(A,B,E)$ be a bipartite graph with proper
edge coloring given by $\chi:E\to\{1,\ldots,d\}$.

\begin{description}
\item[a)]If $G$ does not have a heavy path then
$|E|\le2\sqrt{d|A|\,|B|}\le(|A|+|B|)\sqrt d$.
\item[b)]If $G$ does not have a slow walk then $|E|\le(|A|+|B|)(\log d+2)$.
\item[c)]If $G$ does not have a fast walk then $|E|\le2(|A|+|B|)(\log d+2)$.
\end{description}
\end{lemma}

\proof
For any vertex $z\in A\cup B$ denote by $m(z)=\max(\chi(e))$, where the
maximum is for edges $e$ incident to $z$. For an edge $e=(x,y)\in E$ we define
its {\em weight} to be $w(e)=m(x)-\chi(e)$, while its {\em $B$-weight} is
$w_B(e)=m(y)-\chi(e)$. Clearly, both $w(e)$ and $w_B(e)$ are integers in
$[0,d-1]$.

To prove part a) of the lemma assume $G$ does not contain a heavy path. We
set a threshold parameter $t=\lfloor\sqrt{d|A|/|B|}\rfloor$ and  call an edge
$e$ {\em $B$-light} if $w_B(e)<t$, otherwise $e$ is {\em $B$-heavy}.

All the edges incident to a vertex $y\in B$ have different colors, thus they
also have different $B$-weights, so at most $t$ of them can be $B$-light.
The total number of $B$-light edges is at most $|B|t\le\sqrt{d|A|\,|B|}$.

Now consider two edges $e_1=(x,y_1)$ and $e_2=(x,y_2)$ incident to vertex
$x\in A$. Assume $\chi(e_2)\le \chi(e_1)$. If $w_B(e_2)>0$ we can extend the
path formed by these two edges with the edge $e_3$ incident to
$y_2$ having maximum color $\chi(e_3)=m(y_2)$. Clearly,
$\chi(e_3)=\chi(e_2)+w_B(e_2)$ and the resulting path is heavy unless
$\chi(e_2)>\chi(e_1)+w(e_2)$. Therefore, the number of $B$-heavy edges
incident to $x$ is at most $d/(t+1)$. The total number of $B$-heavy edges is
at most $|A|d/(t+1)\le\sqrt{d|A|\,|B|}$.

For the total number of edges we have
$|E|\le2\sqrt{d|A|\,|B|}\le(|A|+|B|)\sqrt d$.
 
For parts b) and c) of the lemma consider an edge $e=(x,y)\in E$. If
$w(e)=0$ we call the edge $e$ {\em maximal}. Clearly, there are at
most $|A|$ maximal edges. If $e$ is not maximal we define $n(e)$ to be
the ``next larger colored edge at $x$'', i.e., $n(e)$ is the edge in $E$
having minimal color $\chi(n(e))$ among edges incident to $x$ and satisfying
$\chi(n(e))>\chi(e)$. We define the {\em gap} of $e$ to be
$g(e)=\chi(n(e))-\chi(e)$. Clearly, $0<g(e)\le w(e)$. We call the edge $e$
{\em heavy} if $w(e)>2g(e)$, otherwise $e$ is {\em light}. Recall, that for
maximal edges $n(e)$ and $g(e)$ are not defined and maximal edges are
neither light nor heavy.

Let $e_1$ and $e_2$ be distinct edges in $E$ incident to a vertex $x\in A$.
If $\chi(e_1)<\chi(e_2)$ then $w(e_1)\ge w(e_2)+g(e_1)$. If $e_1$ is light,
$w(e_1)\ge2w(e_2)$ follows, therefore at most $\lceil\log d\rceil$ light edges are
incident to $x\in A$. Thus the total number of light edges in $E$ is at most
$\lceil\log d\rceil|A|$.

For part b) of the lemma assume $G$ does not contain a slow walk.
Let $e_2=(x_2,y)$ and $e_3=(x_3,y)$ be distinct non-maximal edges in $E$ and
assume $\chi(e_2)<\chi(e_3)$. Let $e_1=n(e_2)$ and let $e_4$ be the maximal
edge incident to $x_3$ in $G$. We have $\chi(e_1)=\chi(e_2)+g(e_2)$ and
$\chi(e_4)=\chi(e_3)+w(e_3)$. The edges $e_1$, $e_2$, $e_3$, and $e_4$
cannot form a slow walk. As $\chi(e_1)>\chi(e_2)<\chi(e_3)<\chi(e_4)$ we must
have $\chi(e_4)<\chi(e_1)$. This implies
$g(e_2)>w(e_3)$, and if $e_2$ is heavy $w(e_2)>2w(e_3)$. Therefore, at
most  $\lceil\log d\rceil$ heavy edges can be incident to $y\in B$. The total
number of heavy edges in $G$ is at most $(\lceil\log d\rceil|B|$.

For the total number of edges we add the bound obtained for light, heavy, and
maximal edges and get $|E|\le(|A|+|B|)\lceil\log d\rceil+|A|$.

Finally for part c) of the lemma we assume $G$ does not contain a fast walk.
Let $E_1$ consist of the edges $(x,y)\in E$ for which $m(x)\ge m(y)$. Assume
without loss of generality that $|E_1|\ge|E|/2$. If this is not the case
consider the same graph with its sides switched. We use here that the
definition of a fast walk and the claimed bound on the number of edges are
both symmetric in the color classes.

As in the previous case consider two non-maximal edges $e_2=(x_2,y)$
and $e_3=(x_3,y)$ in $E_1$ with $\chi(e_2)<\chi(e_3)$. Let $e_1$ be the
maximal edge incident to $x_2$ and let $e_4=n(e_3)$. As $\chi(e_2)<
\chi(e_3)<\chi(e_4)$ but $G$ does not contain a fast walk we must have
$\chi(e_1)<\chi(e_4)$. As $e_2\in E_1$ we must also have $\chi(e_1)=m(x_2)\ge
m(y)\ge \chi(e_3)$. If $e_3$ is heavy we also have
$$\begin{array}{rcl}
\chi(e_3)+w(e_3)-m(y)&>&\chi(e_3)+2g(e_3)-m(y)\\
&\ge&2(\chi(e_3)+g(e_3)-m(y))\\
&=&2(\chi(e_4)-m(y))\\
&>&2(\chi(e_1)-m(y))\\
&=&2(\chi(e_2)+w(e_2)-m(y)).
\end{array}$$
For all the heavy edges $e\in E_1$ incident to $y\in B$ the values
$\chi(e)+w(e)-m(y)$ increase strictly more than by a factor of $2$. As these
values are integers from $[0,d-1]$, there are at most $\lceil\log d\rceil$ heavy edges
in $E_1$ incident to $y$. The total number of heavy edges in $E_1$ is at most
$\lceil\log d\rceil|B|$.

For the total number of edges in $E_1$ we sum our bound on heavy
edges in $E_1$ and the bounds on the light and maximal edges in
$E$. We obtain $|E_1|\le(|A|+|B|)\lceil\log d\rceil+|A|$. Finally we get
$|E|\le2(|A|+|B|)\lceil\log d\rceil+2|A|$. \qed

\end{document}